\documentclass[12pt,leqno]{article}
\usepackage{amssymb}
\usepackage{amsmath}
\usepackage{url}
 \newtheorem{thm}{Theorem}[section]

 \numberwithin{equation}{section}
\newcommand{\Ll}{\mbox{$\mathcal L$}}
\newcommand{\Th}{\mbox{\sf Th}}
\newcommand{\bbb}{\mbox{\sf big}}
\newcommand{\partition}{\mbox{\sf partition}}
\newcommand{\n}{\mbox{(}}
\newcommand{\z}{\mbox{)}}
\newcommand{\Mm}{\mbox{$\mathfrak M$}}
\newcommand{\mn}{\mbox{\sf mn}}
\newcommand{\id}{\mbox{\sf id}}
\newcommand{\di}{\mbox{\sf di}}
\newcommand{\Rel}{\mbox{\sf Rel}}
\newcommand{\choice}{\mbox{\sf choice}}
\newcommand{\C}{\mbox{\sf C}}

\begin{document}

\title{Finite-variable logics do not have weak Beth definability property}
\author{Andr\'eka, H. and N\'emeti, I.}
\date{September 17, 2014}
\maketitle

\begin{abstract}
We prove that $n$-variable logics do not have the weak Beth
definability property, for all $n\ge 3$. This was known for $n=3$
(Ildik\'o Sain and Andr\'as Simon \cite{Simon}), and for $n\ge 5$
(Ian Hodkinson, \cite{Hodk}). Neither of the previous proofs works
for $n=4$. In this paper we settle the case of $n=4$, and we give a
uniform, simpler proof for all $n\ge 3$. The case for $n=2$ is still
open.
\end{abstract}

\section{Introduction}
Definability theory is one of the most exciting and important parts
of logic. It concerns concept formation and structuring our
knowledge by investigating the category of theories. Implicit
definitions are important in understanding concept formation and
explicit definitions are vital ingredients of interpretations
between theories. This has applications in the methodology of
sciences \cite{comparing, vB82, MadDis}.

Beth definability theorem for first-order logic (FOL) states that
each implicit definition is equivalent to an explicit one, modulo
theories. Investigating whether this theorem holds for fragments of
first-order logic gives information about complexity of the explicit
definition equivalent to the implicit one. Beth definability
property is equivalent to surjectivity of epimorphisms in the
associated class of algebras (a theorem of N\'emeti \cite{Nemeti},
see also \cite{Hoogland, BH06, Sain}).

Failure of Beth definability property for the finite variable
fragments was first proved in 1983 \cite{ACN} (for all $n\ge 2$) by
showing that epimorphisms are not surjective in finite-dimensional
cylindric algebras, see \cite{ACMNS}. That proof, translated to
logic, relies inherently on the fact that the implicit definition it
uses is not satisfiable in each model of the theory. The question
came up whether the so-called weak Beth definability property holds
for finite-variable fragments. Weak Beth definability property
differs from the original Beth definability property in that we
require not only the uniqueness, but also the existence of the
implicitly defined relation. In some sense, the weak Beth
definability property is more intuitive, and is considered to be
more important than the (strong) Beth definability property, see
e.g., \cite{Barwise}.

In this paper we prove that $n$-variable logics do not have the weak
Beth definability property either, for all $n\ge 3$. This means that
there are a first-order logic theory, and an implicit definition
that has exactly one solution in each model of the theory, such that
both the theory and the implicit definition are written up with
using $n$ variables only, yet any explicit definition equivalent to
this implicit one has to use more than $n$ variables. For more on
finite variable logics and the Beth definability properties see
\cite{Hodk} and the remarks at the end of this paper.

\section{The Main Theorem}

The $n$-variable fragment $\Ll_n$ of a FOL language $\Ll$, where $n$
is any finite number, is the set of all formulas in $\Ll$ which use
$n$ variables only (free or bound). To make this more concrete, we
may assume that $\Ll$ uses the variables $v_0, v_1,..., $ while
$\Ll_n$ uses only the variables $v_0,v_1,...,v_{n-1}$. In finite
variable fragments we do not allow function or constant symbols, but
we allow equality. Here is a definition of the formulas of $\Ll_n$:
\begin{description}
\item
$R(v_{i_1},...,v_{i_k})$ is a formula of $\Ll_n$ if $R$ is a
$k$-place relation symbol and $i_1,...,i_k<n$.
\item
$v_i=v_j$ is a formula of $\Ll_n$ if $i,j<n$.
\item
$\neg\varphi$,\ \ $\varphi\land\psi$, \ \ $\exists v_i\,\varphi$ are
formulas of $\Ll_n$  whenever\ \ $\varphi,\psi$ are formulas of
$\Ll_n$ and $i<n$.
\end{description}
\noindent The above are all the formulas of $\Ll_n$. We use other
logical connectives, e.g., $\forall v_i, \lor, \to$ as derived ones.
Models, satisfiability of formulas under evaluations of the
variables, validity in $\Ll_n$ are the same as in FOL. The following
theorem says that $\Ll_n$ does not have even the weak Beth
Definability Property whenever $n\ge 3$.:

\begin{thm}\label{nowBforLn}(No weak Beth Property for $\Ll_n$.)
Let $n\ge 3$. There are a theory $\Th$ in the language of an
$n$-place relation symbol $R$ and a binary relation symbol $S$, and
a theory $\Sigma(D)$ in the language of $\Th$ enriched with a unary
relation symbol $D$ such that
\begin{description}
\item in each model of $\Th$ there is a unique relation $D$ for which
$\Sigma(D)$ holds (we call such $\Sigma(D)$ a strong implicit
definition of $D$ in $\Th$)
\item there is no explicit definition for $D$
in $\Th$, i.e., for each $n$-variable formula $\varphi$ in the
language of $\Th$ we have
\[ \Th\cup\Sigma(D) \not\models \forall v_0[
D(v_0)\leftrightarrow\varphi] \ .\] \end{description}
\end{thm}
\bigskip

\noindent {\sf Proof.} We write out the proof in detail for $n=3$.
Generalizing this proof to all $n\ge 3$ will be easy. We will often
write $x,y,z$ for $v_0,v_1,v_2$ and we will write simply $R$ for
$R(x,y,z)$. We will use $U_0(x), U_1(y), U_2(z)$ to be abbreviations
of the formulas on the right-hand sides of the respective
$\leftrightarrow$'s below:
\[ U_0(x):\leftrightarrow \exists yzR,\qquad
U_1(y):\leftrightarrow \exists xzR,\qquad U_2(z):\leftrightarrow
\exists xyR . \] These formulas express the domain of $R$, i.e., the
first projection of $R$, and the second and third projections of
$R$. We will include formulas into $\Th$ that express that $U_0,
U_1, U_2$ are sets of cardinalities $3,2,2$ respectively, and they
form a partition of the universe. We will formulate these properties
with 3 variables after describing the main part of the construction.
Let us introduce the abbreviations $T$ and $\bbb(R)$ as
\[ \begin{array}{l}
T:\leftrightarrow U_0(x)\land U_1(y)\land U_2(z),\mbox{\ \ and}\\
\bbb(R):\leftrightarrow \bigwedge\{ \exists v_iR\leftrightarrow
\exists v_i(T\land\neg R) : i=0,1,2\} .
\end{array}  \]
In the above, $T$ is the ``rectangular hull" of $R$, and $\bbb(R)$
expresses that $R$ cuts this hull into two parts each of which is
sensitive in the sense that as soon as we quantify over them, the
information on how $R$ cuts $T$ into two parts disappears. (Note
that $\bbb(R)$ implies that $\exists v_iR\leftrightarrow \exists
v_iT\leftrightarrow\exists v_i(T\land\neg R)$.) Assume that
$|U_0|=3, |U_1|=2, |U_2|=2$ and $\partition(U_0, U_1, U_2)$ are
formulas in $\Ll_3$ that express the associated meanings. Then we
define
\[ \Th := \{ |U_0|=3, |U_1|=2, |U_2|=2,\quad \partition(U_0, U_1,
U_2),\quad \bbb(R)\} .\] We will show that $\Th$ has exactly one
model, up to isomorphism. But before doing that, let us turn to
expressing the properties we promised about the $U_i$'s with using
three variables.\bigskip

We will use \emph{Tarski's way of substituting} one variable for the
other. I.e., we introduce the abbreviations
\[ U_1\n x\z:\leftrightarrow \exists y(x=y\land U_1(y)),\qquad
U_2\n x\z:\leftrightarrow \exists z(x=z\land U_2(z)) . \] We now can
express that $U_0, U_1, U_2$ form a partition of the universe:
\[ \forall x(U_0(x)\lor U_1\n x\z\lor U_2\n x\z),\qquad \forall x(U_i(x)\to\neg
U_j(x))\quad\mbox{for}\ \ i\ne j,\ \ i,j<3 .\] For expressing the
sizes of the sets $U_i$ we will use the abbreviations
\[ U_1(z):\leftrightarrow \exists z(z=y\land U_1(y)),\qquad
U_2(y):\leftrightarrow \exists z(y=z\land U_2(z)) . \] Now, for
$i=1,2$ we define the formulas
\[\begin{array}{l}
|U_i|\le 2:\leftrightarrow \neg\exists xyz(x\ne y\land x\ne y\land
y\ne z\land U_i(x)\land U_i(y)\land U_i(z)),\\
|U_i|\ge 2:\leftrightarrow \exists xy(x\ne y\land U_i(x)\land
U_i(y)),\\
|U_i|=2:\leftrightarrow |U_i|\ge 2\land|U_i|\le 2.
\end{array}\]
It remains to express that $U_0$ has exactly three elements. In
$\Ll_n$ with $n\ge 4$ we can express $|U_0|=3$ similarly to the
above, but in $\Ll_3$ we have to use another tool. For expressing in
$\Ll_3$ that $U_0$ has exactly 3 elements, we will use the binary
relation $S$. (This is the sole use of $S$ in $\Th$, for $n\ge 4$ we
can omit $S$ from the language.) We are going to express that $S$ is
a cycle of order 3 on $U_0$. The following formulas express that $S$
is a function on $U_0$ without a fixed point:
\[ \forall x\exists y\, S(x,y),\ \ S(x,y)\land S(x,z)\to y=z,\ \ S(x,y)\to (U_0(x)\land U_0(y)\land
x \ne y) .\] The following formula expresses that $U_0$ consists of
exactly one 3-cycle of $S$:
\[ S(x,y)\leftrightarrow \exists z(S(y,z)\land S(z,x)),\quad S(x,y)\lor
S(y,x)\lor x=y .\] In the above, we used Tarski-style substitution
of variables without mentioning (e.g., $U_0(y)$) and we omitted
universal quantifiers in front of formulas (e.g., we wrote
$S(x,y)\land S(x,z)\to y=z$ in place of $\forall xy(S(x,y)\land
S(x,z)\to y=z)$). This expresses that $U_0$ has exactly 3
elements.\bigskip

We turn to showing that $\Th$ has exactly one model up to
isomorphism. Let $\Mm=\langle M, R, S\rangle\models\Th$. Let $U_i,
T$ be defined as above. Then $M$ is the disjoint union of the
$U_i$'s, and the sizes of the $U_i$'s for $i=0,1,2$ are 3,2,2
respectively. (So $M$ has 7 elements.) Let $U_1=\{ b_0, b_1\}$, let
$c,d$ be the two elements of $U_2$ and let \[ X:=\{ u\in U_0 :
\langle u,b_0,c\rangle\in R\} .\] By $\Mm\models\bbb(R)$ and
$|U_2|=2$ we have that $\langle u,b_0,d\rangle\notin R$ if $u\in X$
and $\langle u,b_0,d\rangle\in R$ if $u\in U_0-X$. Hence \[ U_0-X=\{
u\in U_0 : \langle u,b_0,d\rangle\in R\} .\] Also, by
$\Mm\models\bbb(R)$, $X$ has one, or $X$ has two elements (it cannot
be that $X$ has 0 or 3 elements). If $|X|=1$ then let's use the
notation $c_0=c, c_1=d$, and if $|X|=2$ then let $c_0=d, c_1=c$. Let
us name the elements of $U_0$ as $a_0,a_1,a_2$ such that $X=\{a_0\}$
if $|X|=1$, $X=\{ a_1, a_2\}$ if $|X|=2$ and $S=\{ \langle a_i,
a_j\rangle : j=i+1(mod 3)\mbox{\ and\ }i,j\le 3 \}$. This can be
done by $\Mm\models\Th$. The setting so far determines $R$ by
$\Mm\models\bbb(R)$, as follows. For all $i\le 2,\, j,k\le 1$ we
have $\langle a_i,b_j,c_k\rangle\in R$ if and only if $\langle
a_i,b_{j+1(mod 2)},c_k\rangle\in T-R$ if and only if $\langle
a_i,b_j,c_{k+1(mod 2)}\rangle\in T-R$. This is so by
$\Mm\models\bbb(R)$ and by $|U_i|=2$ for $i=1,2$. From this we have
that
\[ \begin{array}{llll}
R = &\{ \langle u,b_i,c_j\rangle : u=a_0&\mbox{ and }i+j=0(mod
2)\}\ \cup\\
&\{ \langle u,b_i,c_j\rangle : u=a_1\lor u=a_2&\mbox{ and }i+j=1(mod
2)\}\ .
\end{array}\]
We have seen that all models of $\Th$ are isomorphic to each other.
The above also show that there is no automorphism of $\Mm$ that
would move $\{ a_0\}$.\smallskip

We are ready to formulate our implicit definition $\Sigma(D)$. We
design $\Sigma(D)$ so that, by using the above notation, it
specifies $\{ a_0\}$. We will write $D$ in place of $D(x)$.
\[ \begin{array}{lll}
\Sigma(D):=\{ &T\land \neg D\land R &\to \ \ \forall x(T\land\neg
D\to
R)),\\
&T\land\neg D\land \neg R &\to \ \ \forall x(T\land\neg D\to
\neg R)),\\
&D\to U_0(x),\quad |D|=&\!\!\!\!\!1\qquad\} .
\end{array}\]
Then in each model of $\Th$ there is exactly one unary relation $D$
for which $\Sigma(D)$ holds, namely $D$ has to be the unary relation
$\{ a_0\}\subseteq U_0$. Thus $\Sigma(D)$ is a strong implicit
definition of $D$ in $\Th$.\bigskip

It remains to show that $\Sigma$ cannot be made explicit in $\Ll_3$,
i.e., there is no 3-variable formula $\varphi$ in the language of
$\Th$ for which $\Th\cup\Sigma(D)\models D\leftrightarrow\varphi$.
Our plan is to list all the $\Ll_3$-definable relations in the above
model and observe that $\{ a_0\}$, the relation $\Sigma$ defines, is
not among them.
For any $\varphi\in\Ll_3$ define
\[ \mn(\varphi) := \{ \langle a,b,c\rangle :
\Mm\models\varphi[a,b,c]\} .\]
In the above, $\Mm\models\varphi[a,b,c]$ denotes that the formula
$\varphi$ is true in $\Mm$ when the variables $v_0, v_1, v_2$ are
evaluated to $a,b,c$ respectively, and $\mn$ abbreviates
``\emph{meaning}". Let \[ A := \{ \mn(\varphi) : \varphi\in\Ll_3\}
.\] Clearly, $A$ is closed under the set Boolean operations because
\[\begin{array}{l}
\mn(\varphi\land\psi) = \mn(\varphi) \cap \mn(\psi), \\
\mn(\neg\varphi) = M^3 - \mn(\varphi),
\end{array}\]
and so $A$ is closed under intersection and complementation w.r.t.\
$M^3$, the set of all $M$-termed 3-sequences. Since $M$ is finite,
this implies that $A$ is atomic and the elements of $A$ are exactly
the unions of some atoms.
\smallskip

We will list all the atoms of $A$. It is easy to see that the
elements $U_i\times U_j\times U_k$ for $i,j,k\le 2$ are all in $A$
and they form a partition of $M^3$. To list the atoms of $A$, we
will list the atoms below each $U_i\times U_j\times U_k$ by
specifying a partition of each. For $i,j,k\le 2$ let's abbreviate
the sequence $\langle i,j,k\rangle$ by $ijk$.

$U_0\times U_1\times U_2$ is $T$, and the partition of $T$ will be
$\{ R, T-R\}$. For $ijk$ a permutation of $012$, the partition of
$U_i\times U_j\times U_k$, the permuted version of $T$, will be the
correspondingly permuted versions of $R$ and $T-R$. Formally:
Assume $i,j,k$ are all distinct, i.e., they form a permutation of
$0,1,2$. We define
\[\begin{array}{l}
X(ijk,r) := \{ \langle u_i,u_j,u_k\rangle : \langle
u_0,u_1,u_2\rangle\in R\} ,\\
X(ijk,-r) := \{ \langle u_i,u_j,u_k\rangle\in U_i\times U_j\times
U_k : \langle u_0,u_1,u_2\rangle\notin R\} .
\end{array}\]
We note that
\[ X(012,r)=R,\quad \mbox{and}\quad X(012,-r)=T-R .\]
Note that
\[ \mn(R(v_i,v_j,v_k)) = X(ijk,r) ,\]
and the same for $-r$ in place of $r$, so $X(ijk,r), X(ijk,-r)$ are
elements of $A$.

Assume now that $ijk$ is not repetition-free, i.e., $|\{
i,j,k\}|<3$. In these cases the blocks of the partition of
$U_i\times U_j\times U_k$ will be put together from partitions of
$U_m\times U_n$ ($m,n<3$). Recall that $S=\{\langle a_0,a_1\rangle,
\langle a_1,a_2\rangle, \langle a_2,a_0\rangle\}$.  We define
\[ \begin{array}{l}
\overline{S} := \{\langle a,b\rangle : \langle b,a\rangle\in S\} ,\\
\id_i := \{\langle a,a\rangle : a\in U_i\} ,\\
\di_i := \{\langle a,b\rangle : a\ne b,\ \ a,b\in U_i\} .
\end{array}\]
Above, $\id_i, \di_i$ abbreviate ``\emph{identity} on $U_i$", and
``\emph{diversity} on $U_i$", respectively, and $\overline{S}$ is
the inverse of $S$. Since $S$ is a cycle on the three-element set
$U_0$, its inverse $\overline{S}$ is its complement in the diversity
element of $U_0$, so $\{ S, \overline{S}, \id_0\}$ is a partition of
$U_0\times U_0$. Also, $\{ \di_i, \id_i\}$ is a partition of
$U_i\times U_i$ for $i=1,2$.
We are ready to define the ``binary partitions" as follows
\[\begin{array}{l}
\Rel_{00}:=\{ S,\overline{S}, \id_0\},\quad \Rel_{11}:=\{ \di_1,
\id_1\}, \quad \Rel_{22}:=\{ \di_2, \id_2 $\}$, \\
\Rel_{ij}:=\{ U_i\times U_j\}\ \ \mbox{for}\ \  i\ne j.
\end{array}\]
Note that for all $e\in\Rel_{ij}$, $e'\in\Rel_{jk}$ we have $e\circ
e'\in\Rel_{ik}$, where $\circ$ denotes the operation of composing
binary relations.
In general, when $|\{ i,j,k\}|<3$ and $e=\langle
e_0,e_1\rangle\in\Rel_{ij}\times\Rel_{jk}$ we define
\[ \begin{array}{l}
X(ijk,e) := \{ \langle a,b,c\rangle\in U_i\times U_j\times U_k :
 \langle a,b\rangle\in e_0,\ \langle b,c\rangle\in e_1\} .
 \end{array}\]
Notice that we already defined $X(ijk,e)$ for the case when $i,j,k$
are distinct and $e\in\{ r,-r\}$. Let $\choice(e,ijk)$ denote
$e\in\{ r,-r\}$ when $ijk$ is repetition-free, and $e=\langle
e_0,e_1\rangle,\ e_0\in\Rel_{ij},\ e_1\in\Rel_{jk}$ otherwise.
Define
\[\begin{array}{l}
B := \{ X(ijk,e) : i,j,k\le 2, \choice(e,ijk)\} ,\\
C := \{\bigcup Y : Y\subseteq B\} .\end{array} \]

The following notation will be convenient when $\choice(e,ijk)$ and
$ijk$ is not repetition-free.
\[\begin{array}{l}
e_{01}:=e_0,\quad e_{12}:=e_1,\quad e_{02}:=e_0\circ e_1,\\
e_{ij}:=\overline{S}\quad\mbox{when $i>j$ and $e_{ji}=S$},\\
e_{ij}:=e_{ji}\quad\mbox{when $i>j$ and $e_{ji}\ne S$}.
\end{array}\]
The intuitive meaning of $e_{ij}$ is that $\langle a_i,
a_j\rangle\in e_{ij}$ whenever $\langle a_0,a_1\rangle\in e_0$ and
$\langle a_1,a_2\rangle\in e_1$.

We want to prove that $A=C$. We show $A\subseteq C$ by showing
$\mn(\varphi)\in C$ for all $\varphi\in\Ll_3$, by induction on
$\varphi$. Atomic formulas: \[\begin{array}{l}
\mn(R(v_i,v_j,v_k))=X(ijk,r)\quad\mbox{when}\ \ |\{i,j,k\}|=3,\\
\mn(R(v_i,v_j,v_k))=\emptyset\quad\mbox{otherwise},\\
\mn(S(v_i,v_j)) = \bigcup\{ X(n_1n_2n_3,e) : n_i=n_j=0,
e_{n_in_j}=S\} ,\\
\mn(v_i=v_j)=\bigcup\{ X(n_1n_2n_3,e) : n_i=n_j,
e_{n_in_j}\in\{\id_0,\id_1,\id_2\}\} .\end{array}\] Clearly, $M^3\in
C$, and $C$ is closed under complementation with respect to $M^3$
and intersection, because $B$ is finite and its elements form a
partition of $M^3$. Thus,
\[\mn(\neg\varphi)\in C,\quad \mn(\varphi\land\psi)\in
C\qquad\mbox{whenever}\ \ \mn(\varphi),\mn(\psi)\in C .\] To deal
with the existential quantifiers, let us define for arbitrary
$H\subseteq M^3$
\[ \begin{array}{l}
\C_0H := \{ \langle a,b,c\rangle\in M^3 : \langle a',b,c\rangle\in H\ \ \mbox{for some} \ a'\} ,\\
\C_1H := \{ \langle a,b,c\rangle\in M^3 : \langle a,b',c\rangle\in H\ \ \mbox{for some} \ b'\} ,\\
\C_2H := \{ \langle a,b,c\rangle\in M^3 : \langle a,b,c'\rangle\in
H\ \
 \mbox{for some} \ c'\} .
\end{array}\]
Then we have, by the definition of the meaning of the existential
quantifiers, that for all $i\le 2$
\[ \mn(\exists v_i\varphi) = \C_i\mn(\varphi) .\]
Thus, to show that
\[ \mn(\exists v_i\varphi)\in C\quad\mbox{whenever}\ \
\mn(\varphi)\in C \] it is enough to show that $C$ is closed under
$\C_i$,  i.e., $\C_iX\in C$ whenever $X\in C$ (and $i\le 2$). Since
$\C_i$ is additive,  i.e., $\C_i(X\cup Y)=\C_i(X)\cup\C_i(Y)$, it is
enough to show that \[ \C_mX(ijk,e)\in C\quad \mbox{for all}\ \
i,j,k,m\le 2,\ \ \mbox{and good choice $e$ for $ijk$}. \] Assume
$i,j,k$ are distinct and $e\in\{ r,-r\}$. Then by
$\Mm\models\bbb(R)$
\[ \begin{array}{l}\C_0X(ijk,e) = M\times U_j\times U_k,\\
\C_1X(ijk,e) = U_i\times M\times U_k,\\
\C_2X(ijk,e) = U_i\times U_j\times M .\end{array}\] It is easy to
check that $U_i\times U_j\times U_k\in C$ for all $i,j,k$, and hence
$V_0\times V_1\times V_2\in C$ whenever the $V_i$ are unions of
$U_0, U_1, U_2$.
When $i,j,k$ are not all distinct
\[\begin{array}{l}
\C_0X(ijk,e) = M\times e_{12} = \{\langle a,b,c\rangle : \langle
b,c\rangle\in e_{12}\}  =\\
\qquad\ \ \bigcup\{ X(mjk,e') : m\le 2, e'_{12}=e_{12}\},\\
\\
\C_1X(ijk,e) = \{\langle a,b,c\rangle : \langle a,c\rangle\in
e_{02}\} =\\
\qquad\ \ \bigcup\{ X(imk,e') : m\le 2, e'_{02}=e_{02}\},\\
\\
\C_2X(ijk,e) = \bigcup\{ X(ijm,e') : m\le 2, e'_{01}=e_{01}\} .
\end{array}\]
We have seen that $A\subseteq C$.

\medskip To show that $C\subseteq A$ we have to check that each
$X(ijk,e)$ is the meaning of a formula $\varphi\in \Ll_3$ in $\Mm$.
We already did this for $X(ijk,r)$, $i,j,k$ distinct. For $ijk=000$
and $e=\langle S,S\rangle$
\[\begin{array}{l}
X(000,\langle S,S\rangle)=
\mn(U_0(x)\land U_0(y)\land U_0(z)\land S(x,y)\land S(y,z)),
\end{array}\]
where $U_0(x)=\exists yzR,\ \ U_0(y)=\exists x(x=y\land U_0(x)),\ \
U_0(z)=\exists x(x=z\land U_0(x))$ are the abbreviations introduced
before. The other cases are similar, we leave checking them to the
reader.\medskip

Finally, to show that $\mn(D(x))=\{\langle a_0,b,c\rangle : b,c\in
M\}\notin A$, observe that the domain of each element in $B$ either
contains $U_0$ or else is disjoint from it, and therefore the same
holds for their unions. Clearly, this is not true for $\mn(D(x))$.
This shows that $\mn(D)\notin A$, i.e., $D$ cannot be explicitly
defined in $\Mm$. Since $\Mm$ is a model of $\Th$, this means that
$\Sigma(D)$ is not equivalent to any explicit definition that
contains only 3 variables.

To generalize the construction and the proof from $n=3$ to $n\ge 4$
is straightforward. In the general case $M$ has $2n+1$ elements, it
is the disjoint union of sets $U_0,U_1,\dots U_{n-1}$ of sizes
$3,2,\dots,2$ respectively and $R=\{ s\in U_0\times\dots\times
U_{n-1} : (s_0=a_0\land\Sigma\{ a_i : 1\le i<n\}\mbox{ is
even})\,\lor\, (s_0\in\{a_1,a_2\}\land\Sigma\{ a_i : 1\le
i<n\}\mbox{ is odd})\}$. \hfill$\Box$
\bigskip\bigskip

There is a FOL-formula $\varphi(v_0)$ for $\Th$ and $\Sigma(D)$ as
in Thm.\ref{nowBforLn} which explicitly defines $D(v_0)$, since the
Beth definability theorem holds for FOL. The above theorem then
implies that this explicit definition has to use more than $n$
variables. Thus, both the theory and the implicit definition use
only $n$ variables, but any equivalent explicit definition has to
use more than $n$ variables. In our example, $D(v_0)$ can be defined
by using $n+1$ variables. Ian Hodkinson \cite{Hodk}, by using a
construction from \cite{GS96}, proved that for any number $k$ there
are also a theory and a (weak) implicit definition using only $n$
variables such that any explicit definition this implicit definition
is equivalent to has to use more than $n+k$ variables.

Theorem~\ref{nowBforLn} implies (the known fact) that Craig's
Interpolation Theorem does not hold for $n$-variable logic, either,
for $n\ge 3$. This is so because in the standard proof of the Beth's
Definability Theorem in, e.g., \cite[Thm.2.2.22]{CK}, the explicit
definition is constructed from an interpoland. Complexity
investigations for Craig's theorem were done earlier, see, e.g.,
Daniel Mundici \cite{Mu82}.

The proof given here proves more than what
Theorem~\ref{nowBforLn} states. In the proof, $\Th$ and $\Sigma(D)$
are written in the so-called {\it restricted} $n$-variable logic,
and $\Sigma(D)$ is not equivalent to any $n$-variable formula using
even infinitary conjunctions and disjunctions in a finite model of
$\Th$. A formula is called restricted if substitution of variables
is not allowed in it, i.e., it uses relational atomic formulas of
form $R(v_0,...,v_k)$ only (and it does not contain subformulas of
form $R(v_{i0},...,v_{ik})$ where $\langle
i0,...,ik\rangle\ne\langle 0,...,k\rangle$), see \cite[Part II,
sec.4.3]{HMT}. Thus the weak Beth definability property fails for a
wide variety of logics, from the restricted $n$-variable fragment
with finite models only, to $L^n_{\infty,\omega}$.

The variant of $\Ll_n$ in which we allow only models of size $\le
n+1$ has the strong Beth definability property, for all $n$, this is
proved in \cite{ACMNS}. Another variant of $\Ll_n$ that has the
strong Beth definability property is when we allow models of all
sizes but in a model truth is defined by using only a set of
selected (so-called admissible) evaluations of the variables (a
generalized model then is a pair consisting of a model in the usual
sense and this set of admissible evaluations). The so-called Guarded
fragments of $n$-variable logics also have the strong Beth
definability property. For more on this see \cite{AvBNJPhL98, G99,
HMO}.

We note that $\Ll_2$ does not have the strong Beth definability
property (this is proved in \cite{ACMNS}), and we do not know
whether it has the weak one. There are indications that it might
have. If so, $\Ll_2$ would be a natural example of a logic
distinguishing the two Beth definability properties. At present, we
only have artificial examples for this, see Chapter XVIII by
Makowsky, J. in \cite[p.689, item 4.2.2(v)]{Barwise}.\bigskip

\noindent{\bf Acknowledgement} This work was completed with the
support of Hungarian National Grant for Basic Research No T81188.

\bigskip

\noindent Alfred Renyi Institute of Mathematics\\
Hungarian Academy of Sciences\\
P.O. Box 127\\
H-1364 Budapest\\
Hungary\\
\smallskip
andreka.hajnal@renyi.mta.hu\\
nemeti.istvan@renyi.mta.hu

\end{document}